\documentclass[12pt]{amsart}

\usepackage{amsmath}
\usepackage{amssymb}
\usepackage{array}

\newtheorem{theorem}{Theorem}

\theoremstyle{definition}

\begin{document}

\title{Generating function identities for $\zeta(2n+2), \zeta(2n+3)$ via the WZ method}

\author{Kh.~Hessami Pilehrood$^{1}$}

\address{Institute for Studies in Theoretical Physics and Mathematics
(IPM), Tehran, Iran} \curraddr{Mathemetics Department, Faculty of
Science, Shahrekord University, Shahrekord, P.O. Box 115, Iran.}
\email{hessamik@ipm.ir, hessamit@ipm.ir, hessamit@gmail.com}
\thanks{$^1$ This research was in part supported by a grant
from IPM (No. 86110025)}

\author{T.~Hessami Pilehrood$^2$}
\thanks{$^2$ This research was in part supported by a grant
from IPM (No. 86110020)}


\date{}


\begin{abstract}
Using the WZ method we present simpler proofs of Koecher's,
Leshchiner's and Bailey-Borwein-Bradley's identities for
generating functions of the sequences $\{\zeta(2n+2)\}_{n\ge 0},
\{\zeta(2n+3)\}_{n\ge 0}.$ By the same method we give several new
representations for these generating functions yielding faster
convergent series for values of the Riemann zeta function.
\end{abstract}

\maketitle

\section{Introduction}
\label{intro}

The Riemann zeta function is defined by the series
$$
\zeta(s)=\sum_{n=1}^{\infty}\frac{1}{n^s}, \qquad \mbox{for} \quad
{\rm Re}(s)>1.
$$
Ap\'ery's irrationality proof of $\zeta(3)$ and series
acceleration formulae for the first values of the Riemann zeta
function going back to Markov's work \cite{ma}
$$
\zeta(2)=3\sum_{k=1}^{\infty}\frac{1}{k^2\binom{2k}{k}}, \qquad
\zeta(3)=\frac{5}{2}\sum_{k=1}^{\infty}\frac{(-1)^{k-1}}{k^3\binom{2k}{k}},
\qquad
\zeta(4)=\frac{36}{17}\sum_{k=1}^{\infty}\frac{1}{k^4\binom{2k}{k}}
$$
stimulated intensive search of similar formulas for other values
$\zeta(n), n\ge 5.$ Many Ap\'ery-like formulae have been proved
with the help of generating function identities (see \cite{ko,
algr, br, ri, bbb}). M. Koecher \cite{ko} (and independently
Leshchiner \cite{le}) proved that
\begin{equation}
\sum_{k=0}^{\infty}\zeta(2k+3)a^{2k}=\sum_{n=1}^{\infty}\frac{1}{n(n^2-a^2)}=
\frac{1}{2}\sum_{k=1}^{\infty}\frac{(-1)^{k-1}}{k^3\binom{2k}{k}}\,\,\frac{5k^2-a^2}%
{k^2-a^2}\,\prod_{m=1}^{k-1}\left(1-\frac{a^2}{m^2}\right),
\label{eq01}
\end{equation}
 for any $a\in {\mathbb C},$ with $|a|<1.$ For even zeta
values,  Leshchiner \cite{le} (in an expanded form) showed that
(see \cite[(31)]{bbb})
\begin{equation}
\sum_{k=0}^{\infty}\left(1-\frac{1}{2^{k+1}}\right)\zeta(2k+2)a^{2k}=\sum_{n=1}^{\infty}
\frac{(-1)^{n-1}}{n^2-a^2}=\frac{1}{2}\sum_{k=1}^{\infty}\frac{1}{k^2\binom{2k}{k}}\,\,
\frac{3k^2+a^2}{k^2-a^2}\,\prod_{m=1}^{k-1}\left(1-\frac{a^2}{m^2}\right),
\label{eq02}
\end{equation}
 for any complex $a,$ with $|a|<1.$ Recently,
D.~Bailey, J.~Borwein and D.~Bradley \cite{bbb} proved another
formula
\begin{equation}
\sum_{k=0}^{\infty}\zeta(2k+2)a^{2k}=\sum_{n=1}^{\infty}\frac{1}{n^2-a^2}=
3\sum_{k=1}^{\infty}\frac{1}{\binom{2k}{k}(k^2-a^2)}\,\prod_{m=1}^{k-1}
\left(\frac{m^2-4a^2}{m^2-a^2}\right), \label{eq03}
\end{equation}
 for any $a\in
{\mathbb C}, |a|<1.$

In this paper,  we   present simpler proofs of identities
(\ref{eq01})--(\ref{eq03}) using the WZ method. By the same method
we give some new representations for the generating functions
(\ref{eq01}), (\ref{eq03}) yielding faster convergent series for
values of the Riemann zeta function.

We recall \cite{z} that a discrete function $A(n,k)$ is called
hypergeometric or closed form (CF) if the quotients
$$
\frac{A(n+1,k)}{A(n,k)} \qquad\mbox{and}\qquad
\frac{A(n,k+1)}{A(n,k)}
$$
are both rational functions of $n$ and $k.$ A pair of CF functions
$F(n,k)$ and $G(n,k)$ is called a WZ pair if
\begin{equation}
F(n+1,k)-F(n,k)=G(n,k+1)-G(n,k). \label{eq04}
\end{equation}
We need the following summation formulas.

{\bf Proposition 1.} (\cite[Formula 2]{az}) {\it For any WZ pair
$(F,G)$
$$
\sum_{k=0}^{\infty}F(0,k)-\lim_{n\to\infty}\sum_{k=0}^nF(n,k)=
\sum_{n=0}^{\infty}G(n,0)-\lim_{k\to\infty}\sum_{n=0}^kG(n,k),
$$
whenever both sides converge. }

{\bf Proposition 2.}  (\cite[Formula 3]{az})  {\it For any WZ pair
$(F,G)$ we have
$$
\sum_{n=0}^{\infty}G(n,0)=\sum_{n=0}^{\infty}(F(n+1,n)+G(n,n))
-\lim_{n\to\infty}\sum_{k=0}^{n-1}F(n,k),
$$
whenever both sides converge. }

As usual,  $(\lambda)_{\nu}$ is the Pochhammer symbol (or the
shifted factorial) defined by
\begin{equation*}
(\lambda)_{\nu}=\frac{\Gamma(\lambda+\nu)}{\Gamma(\lambda)}
=\begin{cases}
 1,     & \quad \nu=0; \\
\lambda(\lambda+1)\ldots (\lambda+\nu-1), & \quad \nu\in {\mathbb
N}.
\end{cases}
\end{equation*}

\section{Proof of Koecher's identity}

Consider
$$
F(n,k)=\frac{(-1)^nk!(1+a)_n(1-a)_n}{(2n+k+1)!((n+k+1)^2-a^2)}.
$$
Then we have
$$
F(n+1,k)-F(n,k)=G(n,k+1)-G(n,k)
$$
with
$$
G(n,k)=\frac{(-1)^nk!(1+a)_n(1-a)_n(5(n+1)^2-a^2+k^2+4k(n+1))}%
{(2n+k+2)!((n+k+1)^2-a^2)(2n+2)},
$$
i.e., $(F,G)$ is a WZ pair and by Proposition 1, we get
$$
\sum_{k=0}^{\infty}F(0,k)=\sum_{n=0}^{\infty}G(n,0),
$$
or
\begin{equation*}
\begin{split}
\sum_{k=1}^{\infty}\frac{1}{k(k^2-a^2)}&=\sum_{n=0}^{\infty}
\frac{(-1)^n(1+a)_n(1-a)_n(5(n+1)^2-a^2)}{(2n+2)!(2n+2)((n+1)^2-a^2)}
\\ &=
\frac{1}{2}\sum_{n=1}^{\infty}\frac{(-1)^{n-1}(5n^2-a^2)}{n^3\binom{2n}{n}(n^2-a^2)}
\prod_{m=1}^{n-1}\left(1-\frac{a^2}{m^2}\right).   \qed
\end{split}
\end{equation*}

\section{Proof of Leshchiner's identity}

Consider
$$
F(n,k)=\frac{(-1)^kk!(1+a)_n(1-a)_n(n+k+1)}{(2n+k+1)!((n+k+1)^2-a^2)}.
$$
Then we have
$$
F(n+1,k)-F(n,k)=G(n,k+1)-G(n,k)
$$
with
$$
G(n,k)=\frac{(-1)^kk!(1+a)_n(1-a)_n(3(n+1)^2+a^2+k^2+4k(n+1))}%
{2(2n+k+2)!((n+k+1)^2-a^2)}
$$
and by Proposition 1, we get
$$
\sum_{k=0}^{\infty}F(0,k)=\sum_{n=0}^{\infty}G(n,0),
$$
or
\begin{equation*}
\begin{split}
\sum_{k=1}^{\infty}\frac{(-1)^{k-1}}{k(k^2-a^2)}&=\frac{1}{2}\sum_{n=0}^{\infty}
\frac{(1+a)_n(1-a)_n(3(n+1)^2+a^2)}{(2n+2)!((n+1)^2-a^2)} \\ &=
\frac{1}{2}\sum_{n=1}^{\infty}\frac{3n^2+a^2}{n^2\binom{2n}{n}(n^2-a^2)}
\prod_{m=1}^{n-1}\left(1-\frac{a^2}{m^2}\right).   \qed
\end{split}
\end{equation*}

\section{Proof of the  Bailey-Borwein-Bradley identity}

Consider
$$
F(n,k)=\frac{n!^2(1+a)_k(1-a)_k(1+2a)_n(1-2a)_n}{(2n)!(1+a)_{n+k+1}(1-a)_{n+k+1}}.
$$
Then we have
$$
F(n+1,k)-F(n,k)=G(n,k+1)-G(n,k)
$$
with
$$
G(n,k)=\frac{(1+a)_k(1-a)_k(1+2a)_n(1-2a)_nn!(n+1)!(3n+3+2k)}%
{(1+a)_{n+k+1}(1-a)_{n+k+1}(2n+2)!},
$$
and $(F,G)$ is a WZ pair. Then
$$
\sum_{k=0}^{\infty}F(0,k)=\sum_{n=0}^{\infty}G(n,0),
$$
and therefore,
\begin{equation*}
\begin{split}
\sum_{k=1}^{\infty}\frac{1}{(k^2-a^2)}&=3\sum_{n=0}^{\infty}
\frac{(1+2a)_n(1-2a)_n(n+1)!^2}{(1+a)_{n+1}(1-a)_{n+1}(2n+2)!}
\\ &=
3\sum_{n=1}^{\infty}\frac{1}{\binom{2n}{n}(n^2-a^2)}
\prod_{m=1}^{n-1}\left(\frac{m^2-4a^2}{m^2-a^2}\right),
\end{split}
\end{equation*}
as required. \qed

\section{New generating function identities for $\zeta(2n+2)$ and
$\zeta(2n+3)$}

\begin{theorem} \label{t1}
Let $a$ be a complex number not equal to a non-zero integer. Then
\begin{equation}
\sum_{k=1}^{\infty}\frac{1}{k(k^2-a^2)}=\sum_{n=1}^{\infty}
\frac{a^4-a^2(32n^2-10n+1)+2n^2(56n^2-32n+5)}{2n^3\binom{2n}{n}\binom{3n}{n}
((2n-1)^2-a^2)(4n^2-a^2)}\prod_{m=1}^{n-1}\left(\frac{a^2}{m^2}-1\right).
\label{eq05}
\end{equation}
\end{theorem}
Expanding both sides of (\ref{eq05}) in powers of $a^2$ and
comparing coefficients of $a^{2n}$ gives Ap\'ery-like series for
$\zeta(2n+3)$ for every non-negative integer $n$ convergent at the
geometric rate with ratio $1/27.$ In particular, comparing
constant terms recovers Amdeberhan's formula \cite{am} for
$\zeta(3)$
$$
\zeta(3)=\frac{1}{4}\sum_{n=1}^{\infty}(-1)^{n-1}
\frac{56n^2-32n+5}{n^3(2n-1)^2\binom{2n}{n}\binom{3n}{n}}.
$$
Similarly, comparing coefficients of $a^2$ gives
$$
\zeta(5)=\frac{3}{16}\sum_{n=1}^{\infty}\frac{(4n-1)(16n^3-8n^2+4n-1)}%
{(-1)^{n-1}n^5(2n-1)^4\binom{2n}{n}\binom{3n}{n}} +\frac{1}{4}
\sum_{n=1}^{\infty}\frac{(-1)^n(56n^2-32n+5)}{n^3(2n-1)^2\binom{2n}{n}\binom{3n}{n}}
\sum_{k=1}^{n-1}\frac{1}{k^2}.
$$
{\bf Proof.}  Consider
$$
F(n,k)=\frac{(-1)^nn!(2n)!k!(1+a)_k(1-a)_k(1+a)_n(1-a)_n(1+a)_{2n}(1-a)_{2n}}%
{(3n)!(2n+k+1)!(1+a)_{2n+k+1}(1-a)_{2n+k+1}}.
$$
Then application of the WZ algorithm produces  the WZ mate
$$
G(n,k)=\frac{(-1)^nk!n!(2n)!(1+a)_k(1-a)_k(1+a)_n(1-a)_n(1+a)_{2n}(1-a)_{2n}}%
{6(3n+2)!(2n+k+2)!(1+a)_{2n+k+2}(1-a)_{2n+k+2}}q(n,k)
$$
satisfying (\ref{eq04}), with
\begin{equation*}
\begin{split}
&q(n,k)=2(2n+1)(a^4-a^2(32n^2+54n+23)+2(n+1)^2(56n^2+80n+29)) \\
&+k^4(9n+6)+k^3(90n^2+132n+48) +k^2(348n^3+792n^2-15a^2n+594n
 +147 \\ &-9a^2)+k(624n^4+1932n^3
+2214n^2-84a^2n^2-117a^2n+1113n+207-39a^2).
\end{split}
\end{equation*}
By Proposition 1, we have
$$
\sum_{k=0}^{\infty}F(0,k)=\sum_{n=0}^{\infty}G(n,0)
$$
or equivalently,
\begin{equation*}
\begin{split}
&\qquad\qquad\qquad\qquad\sum_{k=1}^{\infty}\frac{1}{k(k^2-a^2)}= \\
&\sum_{n=0}^{\infty}
\frac{(-1)^nn!(1-a)_n(1+a)_n(a^4-a^2(32n^2+54n+23)+2(n+1)^2(56n^2+80n+29))}%
{2(3n+3)!((2n+1)^2-a^2)((2n+2)^2-a^2)},
\end{split}
\end{equation*}
and the theorem follows.   \qed

\begin{theorem} \label{t2}
Let $a$ be a complex number not equal to a non-zero integer. Then
\begin{equation}
\sum_{k=1}^{\infty}\frac{1}{k^2-a^2}=\sum_{n=1}^{\infty}\frac{n^2(21n-8)-a^2(9n-2)}%
{\binom{2n}{n}n(n^2-a^2)(4n^2-a^2)}\prod_{k=1}^{n-1}\left(\frac{k^2-4a^2}{(k+n)^2-a^2}
\right). \label{eq06}
\end{equation}
\end{theorem}
Formula (\ref{eq06}) generates Ap\'ery-like series for
$\zeta(2n+2)$ for every non-negative integer $n$ convergent at the
geometric rate with ratio $1/64.$ In particular, it follows that
\begin{equation}
\zeta(2)=\sum_{n=1}^{\infty}\frac{21n-8}{n^3\binom{2n}{n}^3}
\label{eq065}
\end{equation}
 and
$$
\zeta(4)=\sum_{n=1}^{\infty}\frac{69n-32}{4n^5\binom{2n}{n}^3}
-\sum_{n=1}^{\infty}\frac{21n-8}{n^3\binom{2n}{n}^3}\sum_{k=1}^{n-1}
\left(\frac{4}{k^2}-\frac{1}{(k+n)^2}\right).
$$
Another proof of formula (\ref{eq065}) can be found  in \cite[\S
12]{z}.

 {\bf Proof.} Consider
$$
F(n,k)=\frac{n!^2(1+2a)_n(1-2a)_n(1+a)_{n+k}(1-a)_{n+k}}{(2n)!(1+a)_{2n+k+1}(1-a)_{2n+k+1}}.
$$
Application of the WZ algorithm produces the WZ mate
$$
G(n,k)=\frac{n!^2(1+a)_{n+k}(1-a)_{n+k}(1+2a)_n(1-2a)_n}%
{2(2n+1)!(1+a)_{2n+k+2}(1-a)_{2n+k+2}}q(n,k)
$$
satisfying (\ref{eq04}), with
$$
q(n,k)=(n+1)^2(21n+13)-a^2(9n+7)+2k^3+k^2(13n+11)+k(28n^2+48n+20-2a^2).
$$
By Proposition 1,
$$
\sum_{k=0}^{\infty}F(0,k)=\sum_{n=0}^{\infty}G(n,0),
$$
which implies (\ref{eq06}).   \qed

\begin{theorem} \label{t3}
Let $a$ be a complex number not equal to a non-zero integer. Then
$$
\sum_{k=1}^{\infty}\frac{1}{k(k^2-a^2)}= \frac{1}{4}
\sum_{n=0}^{\infty}\frac{(1+a)_n^2(1-a)_n^2((n+1)^2(30n+19)-a^2(12n+7))}%
{(1+a)_{2n+2}(1-a)_{2n+2}(n+1)(2n+1)}.
$$
\end{theorem}

{\bf Proof.} Consider
$$
F(n,k)=\frac{(1+a)_k(1-a)_k(1+a)_n^2(1-a)_n^2}{(1+a)_{2n+k+1}(1-a)_{2n+k+1}(n+k+1)}.
$$
Then application of the WZ algorithm produces the WZ companion
$$
G(n,k)=\frac{(1+a)_k(1-a)_k(1+a)_n^2(1-a)_n^2q(n,k)}%
{4(1+a)_{2n+k+2}(1-a)_{2n+k+2}(n+k+1)(n+1)(2n+1)},
$$
with
\begin{equation*}
\begin{split}
q(n,k)&=(n+1)^3(30n+19)-a^2(n+1)(12n+7)+2k^3(n+1)+2k^2(7n^2+13n+6)
\\ &+k(34n^3+93n^2+84n-4a^2n+25-3a^2).
\end{split}
\end{equation*}
Now by Proposition 1, the theorem follows.  \qed

\begin{theorem} \label{t4}
Let $a$ be a complex number not equal to a non-zero integer. Then
\begin{equation}
\sum_{k=1}^{\infty}\frac{1}{k(k^2-a^2)}=2\sum_{n=1}^{\infty}
\frac{(-1)^{n-1}}{n^3\binom{2n}{n}^5}\frac{p(n,a)}{(n^2-a^2)(4n^2-a^2)}
\prod_{m=1}^{n-1}\left(\frac{(1-a^2/m^2)^2}{1-a^2/(n+m)^2}\right),
\label{eq07}
\end{equation}
\end{theorem}
where
$$
p(n,a)=a^4-a^2(62n^2-40n+8)+n^2(205n^2-160n+32).
$$
Formula (\ref{eq07}) generates Ap\'ery-like series for
$\zeta(2n+3),$ $n\ge 0,$ convergent at the geometric rate with
ratio $2^{-10}.$ In particular, if $a=0$ we get the formula of
Amdeberhan and Zeilberger \cite{az}
$$
\zeta(3)=\frac{1}{2}\sum_{n=1}^{\infty}\frac{(-1)^{n-1}(205n^2-160n+32)}%
{n^5\binom{2n}{n}^5}.
$$
Comparing coefficients of $a^2$ leads to
\begin{equation*}
\begin{split}
\zeta(5)&=\sum_{n=1}^{\infty}\frac{(-1)^n(31n^2-20n+4)}%
{n^7\binom{2n}{n}^5} \\ &+\sum_{n=1}^{\infty}\frac{(-1)^n
(205n^2-160n+32)}{n^5\binom{2n}{n}^5}\left(\sum_{m=1}^{n-1}
\frac{1}{m^2}-\sum_{m=0}^n \frac{1}{2(m+n)^2}\right).
\end{split}
\end{equation*}

{\bf Proof.} Consider
$$
F(n,k)=\frac{(-1)^k(1+a)_k(1-a)_k(1+a)_n^2(1-a)_n^2(2n-k-1)!k!n!^2}%
{2(n+k+1)!^2(2n)!(1+a)_{2n}(1-a)_{2n}}.
$$
Then
$$
G(n,k)=\frac{(-1)^k(1+a)_k(1-a)_k(1+a)_n^2(1-a)_n^2(2n-k)!k!n!^2q(n,k)}%
{4(2n+1)!(n+k+1)!^2(1+a)_{2n+2}(1-a)_{2n+2}},
$$
with
$$
q(n,k)=(n+1)^3(30n+19)-a^2(n+1)(12n+7)+k(21n^3+55n^2+47n+13-3a^2n-a^2),
$$
is a WZ mate such that
$$
\sum_{n=0}^{\infty}G(n,0)=\sum_{n=0}^{\infty}
\frac{(1+a)_n^2(1-a)_n^2((n+1)^2(30n+19)-a^2(12n+7))}%
{4(n+1)(2n+1)(1+a)_{2n+2}(1-a)_{2n+2}} =\sum_{k=1}^{\infty}
\frac{1}{k(k^2-a^2)},
$$
by Theorem \ref{t3}. Now by Proposition 2, the theorem follows.
\qed

\end{document}